\documentclass[12pt]{article}
\usepackage{url}
\usepackage{hyperref}
\usepackage{amsmath}
\usepackage{amsfonts}
\usepackage{amsthm}
\usepackage{graphicx}
\usepackage{epstopdf}
\usepackage{cite}
\usepackage[top=1in, bottom=1in, left=1in, right=1in]{geometry}
\title{The Linearly Independent Non Orthogonal yet Energy Preserving (LINOEP) vectors}
\author{Pushpendra Singh$^{1,3,}$\footnote{Corresponding author's E-mail address: \texttt{pushpendra.singh@ee.iitd.ernet.in; spushp@gmail.com} (P. Singh).  $^1$\url{http://ee.iitd.ernet.in/} ; $^2$\url{http://www.samsung.com/in/} ; $^3$\url{www.jiit.ac.in/}} , S D Joshi$^{1}$, R K Patney$^{1}$ and Kaushik Saha$^{2}$
\\{\normalsize $^{1}$Department of EE, Indian Institute of Technology Delhi, India}\\
{\normalsize $^{2}$Samsung R \& D Institute India - Delhi, India}\\
{\normalsize $^{3}$Jaypee Institute of Information Technology - NOIDA, India}}
\providecommand{\keywords}[1]{\textbf{\textit{Index terms---}} #1}
\date{\today}
\begin{document}
\maketitle
%\footnote{An example footnote.}
\begin{abstract}
It is well known that, in any inner product space, a set of linearly independent (LI) vectors can be transformed to a set of orthogonal vectors, spanning the same space, by the Gram-Schmidt Orthogonalization Method (GSOM). In this paper, we propose a transformation from a set of LI vectors to a set of LI non orthogonal yet energy (square of the norm) preserving (LINOEP) vectors in an inner product space and we refer it as LINOEP method. We also show that there are various solutions to preserve the square of the norm.
\end{abstract}
\keywords{The Gram-Schmidt Orthogonalization Method (GSOM), Linearly Independent Non Orthogonal yet Energy Preserving (LINOEP) vectors, Empirical mode decomposition (EMD).}
 \section{Introduction}
The Gram-Schmidt Orthogonalization Method (GSOM) is a process for obtaining a set of orthogonal vectors, from a set of linearly independent (LI) vectors in an inner product space. The old set of LI vectors and the new set of orthogonal vectors have the same linear span. Let $Y = \{\mathbf{y}_1,\mathbf{y}_2,...,\mathbf{y}_n\}$ be a set of $n$ LI vectors. A set of orthogonal vectors $S = \{\mathbf{s}_1,\mathbf{s}_2,...,\mathbf{s}_n\}$ is generated from the set $Y$ as follows (for $k=1,2,\dots,n$): %$\mathbf{s}_k=\mathbf{y}_k-\sum_{i=1}^{k-1}c_{ki}\mathbf{s}_i \Leftrightarrow$
\begin{multline}
\mathbf{s}_k=\mathbf{y}_k-\sum_{i=1}^{k-1}c_{ki}\mathbf{s}_i \Leftrightarrow
 \left[ \begin{array}{c} \mathbf{y}_1 \\ \mathbf{y}_2 \\ \vdots \\ \mathbf{y}_n \end{array} \right] = \begin{bmatrix} 1 & 0 & \dots & 0 \\ c_{21} & 1 & \dots & 0 \\ \vdots & \vdots & \ddots & \vdots \\ c_{n1} & c_{n2} & \dots & 1 \end{bmatrix} \left[ \begin{array}{c} \mathbf{s}_1 \\ \mathbf{s}_2 \\ \vdots \\ \mathbf{s}_n \end{array} \right]\label{gso1}\end{multline}
The $c_{ki}$ is obtained by using inner product $ \langle \mathbf{s}_k,\mathbf{s}_i \rangle =0 \text{, for } k \neq i$, i.e.
%$c_{ki}=\int_{0}^T\mathbf{y}_ks_idt/\int_{0}^Ts_i^2 dt$
$c_{ki} = \frac{ \langle \mathbf{y}_k,\mathbf{s}_i \rangle}{\langle \mathbf{s}_i,\mathbf{s}_i \rangle}$
for $i=1,2,\dots,n \text{, and } k\ge i$.
%where $T$ is the total observation period of the signals.
 By taking sum of all the $n$ equations of \eqref{gso1} along with some simple algebraic manipulations, it can be shown that
\begin{equation}
\sum_{i=1}^{n}\mathbf{y}_i=\sum_{i=1}^{n}c_i{\mathbf{s}}_{i} \label{vf1}
\end{equation}
where $c_i=\sum_{k=i}^{n}c_{ki}$ is sum of $i^{th}$ column of the coefficient matrix of \eqref{gso1}.
 It can be easily shown that $c_{ki}=1$, if $k=i$. From \eqref{vf1}, it is easy to show (Plancherel/Parseval Equality)
\begin{equation}
\left\lVert \sum_{i=1}^{n}\mathbf{y}_i \right\rVert^2= \left\lVert \sum_{i=1}^{n}c_i{\mathbf{s}}_{i} \right\rVert^2 =\sum_{i=1}^{n}|c_i|^2\lVert {\mathbf{s}}_{i}\rVert^2
\end{equation}
From set $Y$, there are $n$ choices for selecting first vector, $n-1$ choices for second vector, $n-2$ choices for third vector and 1 choice for last vector, that means there are $n!$ permutations of the set $Y$, and the GSOM would produce $n!$ orthogonal sets of vectors from a set of $n$ LI vectors.

The empirical mode decomposition (EMD) is an adaptive signal analysis algorithm, introduced in \cite{rs1}, for the analysis of  nonlinear and non stationary  time series. The various variants of EMD algorithm are proposed in literature \cite{sp4,sp5,sp6} and the orthogonal property of intrinsic mode functions (IMFs) are discussed in \cite{sp2,sp3}. The LI non orthogonal yet energy (square of the norm) preserving (LINOEP) class of vectors and the following theorem are proposed, in \cite{sp1}, for the development of energy preserving EMD (EPEMD) algorithm.
\newtheorem{name}{Theorem}
\begin{name}
Let $H$ be a Hilbert space over the field of complex numbers, and let $\{\mathbf{x},\mathbf{x}_1,\cdots,\mathbf{x}_{n}\}$ be a set of vectors satisfying the following conditions:
%\\$(i)$ $\{\mathbf{x}_1,\cdots,\mathbf{x}_{n+1}\}$ are lineally independent (LI).
\begin{equation}
(i) \qquad \qquad \qquad \qquad  {\mathbf{x}_i \perp \sum_{j=i+1}^{n}\mathbf{x}_j} \label{st1}
\end{equation}
\begin{equation}
(ii) \qquad \qquad \qquad \qquad \qquad  \mathbf{x} =\sum_{i=1}^{n}\mathbf{x}_i \label{st2}
\end{equation}
Then in the representation, given in \eqref{st2}, the square of the norm, and hence energy is preserved, i.e.
\begin{equation}
\left\lVert \mathbf{x} \right\rVert^2= \left\lVert \sum_{i=1}^{n}\mathbf{x}_i \right\rVert^2 =\sum_{i=1}^{n}\lVert \mathbf{x}_i\rVert^2
\end{equation}
\end{name}
In this result, pairwise orthogonality is not required and only last two vectors (i.e. $\mathbf{x}_{n-1}$ and $\mathbf{x}_n$) are orthogonal.
It is the Pythagoras's theorem for $n=2$, and it is the Parseval's theorem when all the basis vectors $\mathbf{x}_i$ are orthogonal.
We use the LINOEP class of vectors, the above theorem and the GSOM for development of LINOEP method in next section.
\section{The LINOEP Method}
There is the GSOM for transformation of a set of linearly independent (LI) vectors to a set of orthogonal vectors in an inner product space. We, motivated by the GSOM, propose the transformation from a set of LI vectors to a set of LINOEP vectors in an inner product space.

Let $Y = \{\mathbf{y}_1,\mathbf{y}_2,...,\mathbf{y}_n\}$ be a set of $n$ LI vectors. This method generates a set of LI non orthogonal yet energy preserving (LINOEP) vectors $S = \{\mathbf{c}_1,\mathbf{c}_2,...,\mathbf{c}_n\}$ from a set $Y$ as follows (for $k=1,2,\dots,n-1$): %$c_k=\mathbf{y}_k-\sum_{i=k+1}^{n}\alpha_{k}c_i \Leftrightarrow$
\begin{multline}
\mathbf{c}_k=\mathbf{y}_k-\alpha_{k} \sum_{i=k+1}^{n}\mathbf{c}_i \text{, } \mathbf{c}_n=\mathbf{y}_n  \Leftrightarrow
 \left[ \begin{array}{c} \mathbf{y}_1 \\ \mathbf{y}_2 \\ \vdots \\ \mathbf{y}_n \end{array} \right] = \begin{bmatrix} 1 & \alpha_1 & \dots & \alpha_1 \\ 0 & 1 & \dots & \alpha_2 \\ \vdots & \vdots & \ddots & \vdots \\ 0 & 0 & \dots & 1 \end{bmatrix} \left[ \begin{array}{c} \mathbf{c}_1 \\ \mathbf{c}_2 \\ \vdots \\ \mathbf{c}_n \end{array} \right]\label{li2epb1}\end{multline}
The values of $\alpha_{k}$ are obtained by equation:
\begin{equation}
%\alpha_{i}=\frac{\int_{0}^Ty_i\left[\sum_{j=i+1}^{n} c_j\right]dt}{\int_{0}^T\left[\sum_{j=i+1}^{n} c_j\right]^2 dt} \label{li2epb2}
\alpha_{k}=\frac{\left\langle \mathbf{y}_k, \sum_{i=k+1}^{n} \mathbf{c}_i \right\rangle}{\left\langle \sum_{i=k+1}^{n} \mathbf{c}_i, \sum_{i=k+1}^{n} \mathbf{c}_i \right\rangle} \qquad k=n-1,\dots,2,1 \label{li2epb2}
\end{equation}
where ${\mathbf{c}_i \perp \sum_{j=i+1}^{n}\mathbf{c}_{j}}$, i.e., inner product $\langle {\mathbf{c}_i, \sum_{j=i+1}^{n}\mathbf{c}_{j}} \rangle =0$ for $i=1,2,\dots,n-1$.
%and $T$ is the total observation period of the vector.
The generated vectors $\mathbf{c}_i$ are LINOEP, thus,
\begin{equation}
\left \lVert \sum_{i=1}^{n}{\mathbf{c}}_{i} \right \rVert ^2=\sum_{i=1}^{n}\lVert {\mathbf{c}}_{i}\rVert ^2. \label{li2epb3}
\end{equation}
By taking sum of all the $n$ equations of \eqref{li2epb1}, along with algebraic manipulation, it can be easily shown that
$\sum_{i=1}^{n}\mathbf{y}_i=\sum_{i=1}^{n}{\mathbf{c}}_{i} +\sum_{i=1}^{n-1}\beta_i{\mathbf{c}}_{i+1}$, where $\beta_i=\sum_{j=1}^{i}\alpha_{j}$, for $i=1,\cdots,n-1$.
%is sum of $i^{th}$ column of coefficient matrix of \eqref{li2epb1}.
Let, $\mathbf{p}_2=\sum_{i=1}^{n}{\mathbf{c}}_{i}$ and $\mathbf{p}_1=\sum_{i=1}^{n-1}\beta_i{\mathbf{c}}_{i+1}$, and further let, $\mathbf{z}_1=\mathbf{p}_2$, and $\mathbf{z}_2=\mathbf{p}_1 - \gamma \mathbf{z}_1 $, where $\gamma= \frac{\left\langle \mathbf{p}_1, \mathbf{z}_1 \right\rangle}{\left\langle \mathbf{z}_1, \mathbf{z}_1 \right\rangle}$ such that $\mathbf{z}_1$ and $\mathbf{z}_2$ are orthogonal. Through the addition of $\mathbf{z}_1$ and $\mathbf{z}_2$, we obtain $\mathbf{p}_1+\mathbf{p}_2=(1+\gamma)\mathbf{z}_1+\mathbf{z}_2$, and, $\sum_{i=1}^{n}\mathbf{y}_i=(1+\gamma)\mathbf{z}_1 + \mathbf{z}_2$. It can be easily shown that
\begin{equation}
\sum_{i=1}^{n}\mathbf{y}_i=\sum_{i=1}^{n+1}\mathbf{d}_i, \label{li2epb4}
\end{equation}
where $\mathbf{d}_i=(1+\gamma)\mathbf{c}_i$ for $i=1,\cdots,n$ and $\mathbf{d}_{n+1}=\mathbf{z}_2$. From \eqref{li2epb3},\eqref{li2epb4}, we obtain
\begin{equation}
\left \lVert \sum_{i=1}^{n}{\mathbf{y}}_{i} \right \rVert ^2=\left \lVert \sum_{i=1}^{n+1}{\mathbf{d}}_{i} \right \rVert ^2=\sum_{i=1}^{n+1}\lVert {\mathbf{d}}_{i}\rVert ^2 \label{li2epb5}
\end{equation}
There are $n!$ permutations of the set $Y$, and hence this method would produce $n!$ sets of LINOEP vectors from a set of $n$ LI vectors.
%In the GSOM $[n(n-1)/2]$ coefficients are calculated and linear combination of orthogonal vectors are energy preserving.
 In this method, we calculate $(n)$ coefficients (i.e., $n-1$ values of $\alpha_i$ and one value of $\gamma$) and the only $\mathbf{x}=\delta_1\mathbf{c}_1+\delta_2[\mathbf{c}_2+\mathbf{c}_3+\dots+\mathbf{c}_n]$ linear combinations are energy preserving and all other linear combinations are not energy preserving.

We observe that if vectors $\mathbf{y}_i \in \mathbb{C}^m$ and $m \ge n$, and all the vectors $\{\mathbf{y}_i\}_{i=1}^{n}$ are LI, then ($n$) out of $(n+1)$ generated vectors $\{\mathbf{d}_i\}_{i=1}^{n}$ are LINOEP, and last vector $\mathbf{d}_{n+1}$ is linear combination of $\{\mathbf{d}_i\}_{i=1}^{n}$ vectors, and the complete set is non orthogonal yet energy preserving  (NOEP).
%The following observations have been made:
%\\(1) If vectors $\mathbf{y}_i \in \mathbb{R}^m$ and $m \ge n$, and $(n)$ vectors $\mathbf{y}_i$ are LI, then ($n$) of $(n+1)$ generated vectors $\{\mathbf{d}_i\}_{i=1}^{n}$ are LINOEP, and last vector $\mathbf{d}_{n+1}$ is linear combination of $\{\mathbf{d}_i\}_{i=1}^{n}$ vectors, and the complete set is non orthogonal yet energy preserving  (NOEP).
%\\(2) If vectors $\mathbf{y}_i \in \mathbb{R}^m$ and $m < n$, and $(m)$ vectors $\mathbf{y}_i$ are LI, then $(m)$ out of $(n+1)$ generated vectors $\mathbf{d}_i$ %are LI, and the complete set is NOEP.
%\section{How may solutions are to preserve the energy}
\section{Various solutions to preserve the square of norm}
%From \eqref{li2epb3}, we observe that
\begin{figure}[!b]
\centering
\includegraphics[angle=0,width=0.99\textwidth,height=0.4\textwidth]{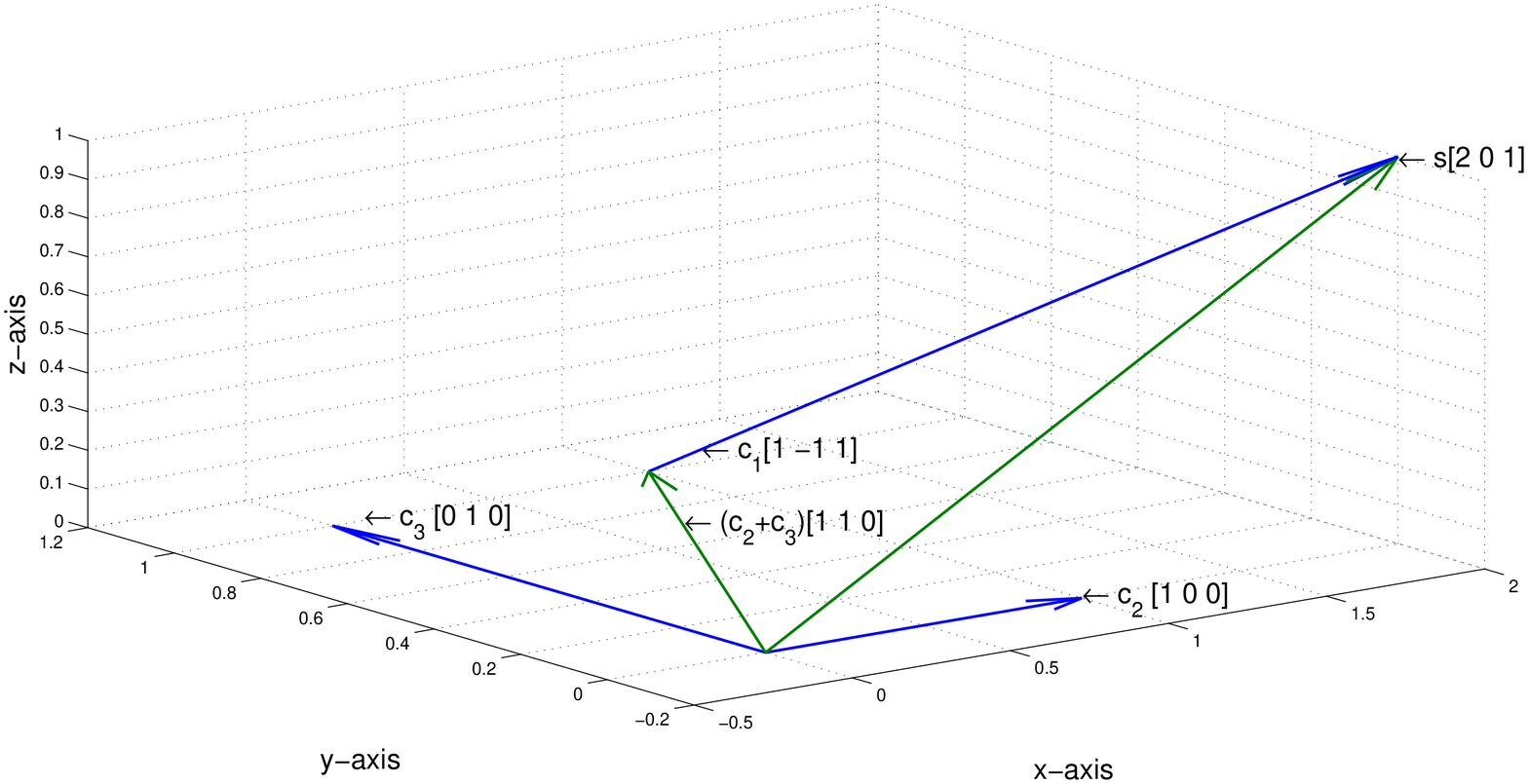}
\caption{Three non orthogonal vectors $c_1, c_2, c_3$ such that $c_1\perp (c_2+c_3)$ and vector $s = c_1 + c_2 + c_3$ in 3-D.}
\label{fig:ex1}
\end{figure}
\begin{figure}[!b]
\centering
\includegraphics[angle=0,width=0.99\textwidth,height=0.4\textwidth]{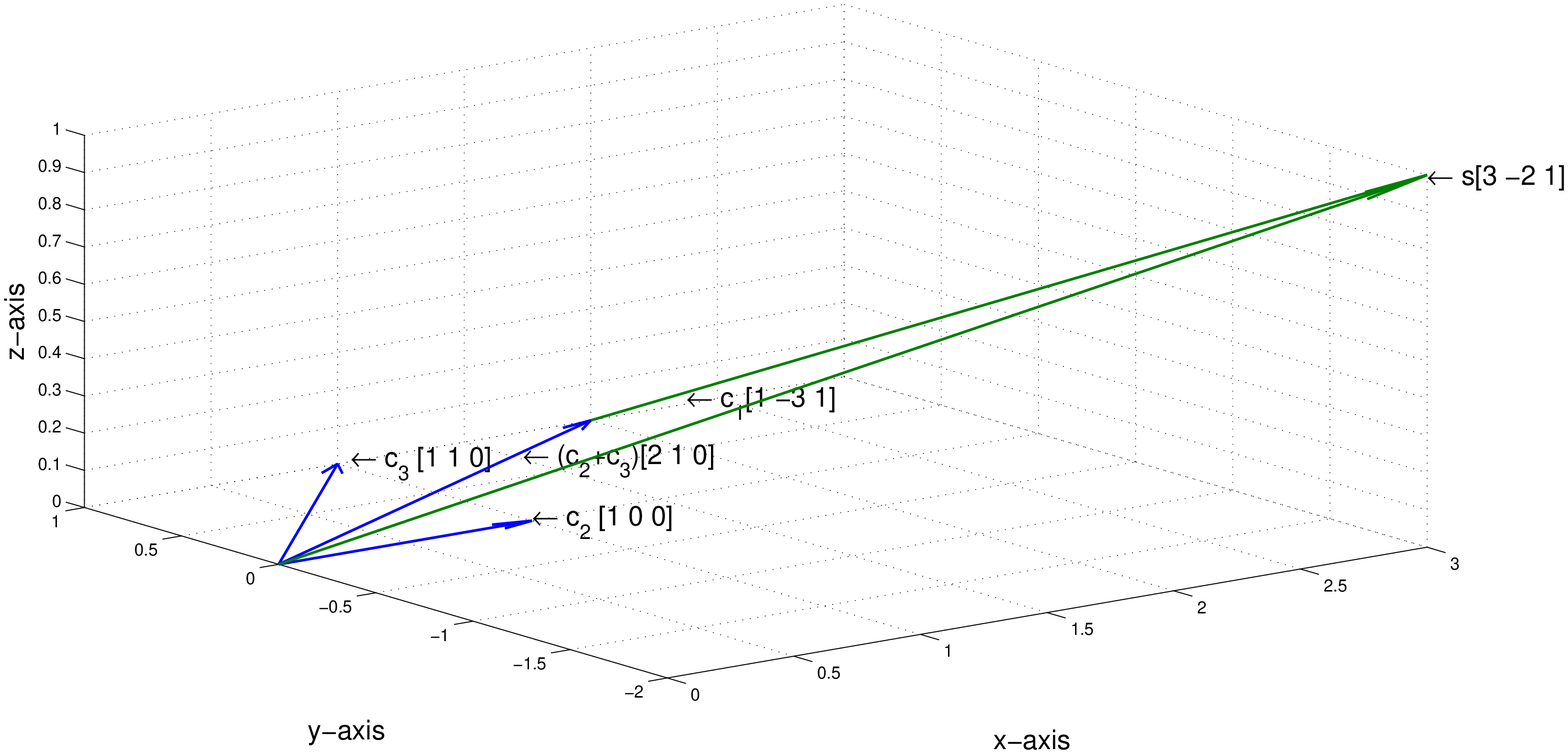}
\caption{Three non orthogonal vectors $c_1, c_2, c_3$ such that $\left\langle \mathbf{c}_{2},\mathbf{c}_{3}\right\rangle+\left\langle \mathbf{c}_{1},\mathbf{c}_{2}\right\rangle=-\left\langle \mathbf{c}_{1},\mathbf{c}_{3}\right\rangle$ and vector $s = c_1 + c_2 + c_3$ in 3-D.}
\label{fig:ex11}
\end{figure}
It is easy to show that
\begin{equation}
\left \lVert \sum_{i=1}^{n}{\mathbf{c}}_{i} \right \rVert ^2=\sum_{i=1}^{n}\lVert {\mathbf{c}}_{i}\rVert ^2 + \sum_{i=1}^{n}\sum_{\substack{j=1\\j\ne i}}^{n}\left\langle {\mathbf{c}}_{i},{\mathbf{c}}_{j}\right\rangle\label{li2epb6}
\end{equation}
and to preserve the energy (square of the norm), cross terms are set to zero, i.e.,
 \begin{equation}
 \sum_{i=1}^{n}\sum_{\substack{j=1\\j\ne i}}^{n}\left\langle {\mathbf{c}}_{i},{\mathbf{c}}_{j}\right\rangle=0. \label{li2epb61}
 \end{equation}
Assuming the underlying field to be $\mathbb{R}$, we write
\begin{equation}
\sum_{i=1}^{n}\sum_{\substack{j=1\\j\ne i}}^{n}\left\langle {\mathbf{c}}_{i},{\mathbf{c}}_{j}\right\rangle=2\sum_{i=1}^{n-1}\sum_{j=i+1}^{n}\left\langle {\mathbf{c}}_{i},{\mathbf{c}}_{j}\right\rangle=0. \label{li2epb7}
\end{equation}
There are various possible solution to make the cross term zero in Eq.~\eqref{li2epb7}:
\\(1) vectors are pairwise orthogonal, i.e., $\left\langle {\mathbf{c}}_{i},{\mathbf{c}}_{j}\right\rangle=0$ for $i\ne j$.
\\(2) use Eq.~\eqref{st1}, i.e, $\sum_{i=1}^{n-1}\sum_{j=i+1}^{n}\left\langle {\mathbf{c}}_{i},{\mathbf{c}}_{j}\right\rangle=0$. There are $n$ choices for selecting first vector, $n-1$ choices for second vector, $n-2$ choices for third vector and $3$ choices for last three vectors, hence, this gives $\frac{n!}{2}$ ways to make cross term zero.
\\(3) There are so many other ways to make the cross term zero which can be obtained.
%\\(3) There are so many other ways which is open question. How many possible ways are to make the corss term zero and how to construct them?
\\\textbf{Example 1.} For $n=2$, there is only one solution to make the corss term zero and this is $\left\langle {\mathbf{c}}_{1},{\mathbf{c}}_{2}\right\rangle=0$.

\noindent \textbf{Example 2.} For $n=3$, there are 7 possible solutions to make the corss term zero:
\\(1) Obvious solution is $\left\langle {\mathbf{c}}_{1},{\mathbf{c}}_{2}\right\rangle=0$, $\left\langle {\mathbf{c}}_{2},{\mathbf{c}}_{3}\right\rangle=0$, $\left\langle {\mathbf{c}}_{1},{\mathbf{c}}_{3}\right\rangle=0$, and this gives one solution.
\\(2) use Eq.~\eqref{st1} and this gives three solutions: (a) $\left\langle \mathbf{c}_{1},\mathbf{c}_{2}+\mathbf{c}_{3}\right\rangle=0$ and $\left\langle \mathbf{c}_{2},\mathbf{c}_{3}\right\rangle=0$ (e.g. see Figure~\ref{fig:ex1}). (b) $\left\langle \mathbf{c}_{2},\mathbf{c}_{1}+\mathbf{c}_{3}\right\rangle=0$ and $\left\langle \mathbf{c}_{1},\mathbf{c}_{3}\right\rangle=0$. (c) $\left\langle \mathbf{c}_{3},\mathbf{c}_{1}+\mathbf{c}_{2}\right\rangle=0$ and $\left\langle \mathbf{c}_{1},\mathbf{c}_{2}\right\rangle=0$.
\\(3) obtain vectors such that inner product of two vectors cancel all others, i.e., $\left\langle \mathbf{c}_{1},\mathbf{c}_{2}\right\rangle+\left\langle \mathbf{c}_{1},\mathbf{c}_{3}\right\rangle+\left\langle \mathbf{c}_{2},\mathbf{c}_{3}\right\rangle=0$, this gives three solutions: (a) $\left\langle \mathbf{c}_{1},\mathbf{c}_{2}\right\rangle+\left\langle \mathbf{c}_{1},\mathbf{c}_{3}\right\rangle=-\left\langle \mathbf{c}_{2},\mathbf{c}_{3}\right\rangle$. (b) $\left\langle \mathbf{c}_{1},\mathbf{c}_{3}\right\rangle+\left\langle \mathbf{c}_{2},\mathbf{c}_{3}\right\rangle=-\left\langle \mathbf{c}_{1},\mathbf{c}_{2}\right\rangle$. (c) $\left\langle \mathbf{c}_{2},\mathbf{c}_{3}\right\rangle+\left\langle \mathbf{c}_{1},\mathbf{c}_{2}\right\rangle=-\left\langle \mathbf{c}_{1},\mathbf{c}_{3}\right\rangle$ (e.g. see Figure~\ref{fig:ex11}).

Similarly, when we assume the underlying field to be $\mathbb{C}$, various possible solutions can be obtaied to make the cross term zero in Eq.~\eqref{li2epb61}.
\section{Conclusions}
We have proposed the transformation from a set of linearly independent (LI) vectors to a set of LI non orthogonal yet energy (square of the norm) preserving (LINOEP) and non orthogonal yet energy preserving (NOEP) vectors in an inner product space. We have also shown that there are various solutions to preserve the square of the norm.
\section*{Acknowledgment}
We would like to thank JIIT Noida, for permitting to carry out research at IIT Delhi and providing all required resources throughout this study.

\end{document}